\documentclass[a4paper,leqno,12pt,twoside]{amsart}
\usepackage{amsmath,amscd,amssymb,upgreek} 
\usepackage[all]{xy}
\usepackage{graphicx}
\newtheorem{Lem}{Lemma \thesection.\!\!}
\newtheorem{Th}[Lem]{Theorem \thesection.\!\!}

\newtheorem{Def}[Lem]{Definition \thesection.\!\!}
\newtheorem{Ex}[Lem]{Example \thesection.\!\!}

\newtheorem{Prop}[Lem]{Proposition \thesection.\!\!}
\newtheorem{Rem}[Lem]{Remark \thesection.\!\!}
\newtheorem{Lem and Def}[Lem]{Lemma and Definition \thesection.\!\!}

\newtheorem{Prop and Def}[Lem]{Proposition and Definition   \thesection.\!\!}

\def\a{\alpha }
\def\b{\beta }

\def\B{\mathcal{B}}

\def\G{\Gamma }

\def\D{\mathcal D }
\def\e{\epsilon }

\def\l{\lambda }

\def\S{\mathcal S }
\def\s{\sigma }

\def\z{\zeta}

\def\ff{\varphi }
\def\f{\mathop{\varphi}}

\newcommand{\A}{\mathcal{A}}
\newcommand{\C}{\mathbb{C}}
\newcommand{\N}{\mathbb{N}}

\newcommand{\Z}{\mathbb{Z}}

\begin{document}
\title{Logarithmic Moduli Spaces for Surfaces of Class $\rm VII$}
\author{Karl OELJEKLAUS and Matei TOMA}  
\date{\today}

\thanks{Part of this work was done while the first author visited 
the University of Osnabr\"uck under the program  "Globale Methoden in der komplexen
Geometrie" of the DFG and while the second author visited the
Max-Planck-Institut f\"ur Mathematik in Bonn and the LATP, Universit\'e de
Provence. We thank these institutions for their hospitality and for financial support. \\
Furthermore the authors wish to thank Georges Dloussky for numerous discussions on surfaces of class $\rm VII$}

\address{{\it Karl Oeljeklaus}:
LATP-UMR(CNRS) 6632, 
CMI-Universit\'e d'Aix-Marseille I, 39, rue Joliot-Curie, F-13453 Marseille Cedex 13, France.}
\email{karloelj@cmi.univ-mrs.fr}
\address{{\it Matei Toma}: Institut f\"ur Mathematik, Universit\"at Osnabr\"uck,
49069 Osnabr\"uck, Germany and 
Institute of Mathematics of the Romanian Academy.}
\email{matei@mathematik.Uni-Osnabrueck.DE}
\urladdr{http://www.mathematik.uni-osnabrueck.de/staff/phpages/tomam.rdf.html}

\begin{abstract}
  In this paper  we describe logarithmic moduli spaces of pairs $(S,D)$ 
consisting of a minimal surface $S$ of class $VII$ with second Betti number $b_2>0$
together with a reduced maximal divisor $D$ of $b_2$ rational curves. The special 
case of Enoki surfaces has already been considered by Dloussky and Kohler.
We use normal forms for the action of the fundamental group of $S\setminus D$ and for 
the associated  holomorphic contraction $(\C,0)\to (\C,0)$.
\end{abstract}

\setcounter{section}{0}
\noindent
\maketitle

\section{Introduction}
Compact complex surfaces with first Betti number $b_1 =1$ 
form the class $\rm VII$ in Kodaira's classification, see \cite{Kod66}.
Minimal such surfaces are said to belong to class $\rm VII_0$. 
When their  second Betti number vanishes 
these surfaces have been completely classified,
see \cite{Inoue74}, \cite{Bog76}, \cite{LYZ94} and \cite{Tel94}. 
Among them 
one finds the Hopf surfaces for which moduli spaces have been constructed
in \cite{Dab82}.\\
In this paper 
we shall restrict our attention to minimal compact complex surfaces
with $b_1=1$ and $b_2>0$. All known surfaces in this subclass contain 
{\it global spherical shells} and can be obtained by a construction due
to Ma. Kato, see \cite{Kato78}. 
A finer subdivision of these surfaces may be done 
by looking at the dual graph 
of the maximal reduced divisor $D$ of rational curves.
One gets:
\begin{enumerate}
\item {\it Enoki surfaces}, 
i.e. the graph of $D$ is a cycle and $D$ is homologically trivial,
\item {\it Inoue-Hirzebruch surfaces}, i.e. the graph of $D$ consists of
one or two cycles and $D$ is not homologically trivial,
(in fact $D$ is an exceptional divisor),
\item {\it intermediate surfaces}, i.e. the graph consists of a cycle with at least
one non-empty tree appended, (in this case too, $D$ is exceptional).
\end{enumerate}

Our point of interest is to study {\it logarithmic moduli spaces}, i.e.
such that the
 maximal reduced divisor $D$ is preserved.
Logarithmic moduli spaces of Enoki surfaces were described in \cite{DK98},
whereas Inoue-Hirzebruch surfaces are logarithmically rigid.\\

In this paper we introduce logarithmic moduli spaces for
intermediate surfaces. This will be done by examining normal
forms of germs of analytic mappings which can be associated to
these surfaces and of the fundamental group of the complement
of the maximal reduced divisor, respectively.\\

The article is organized as follows. General facts on surfaces with global spherical shells are recalled in Section  \ref{facts}. 
Then we give a description of the fundamental group of the complement of the rational curves of an intermediate surface. This description will allow us a better understanding of the isomorphisms between the different parameter spaces we will get for logarithmic deformations.
In Section \ref{germs} we introduce the main ingredients of our proof: the normal form of a contracting germ associated to an intermediate surface and its type. The next two sections are devoted to the decomposition of such normal forms and to the associated blowing up sequences which are needed in order to recover the surface from its contracting germ in a canonical way.
Using these techniques we come in Section \ref{secmoduli} to the description of logarithmically versal families and moduli spaces of intermediate surfaces, see Theorems 7.\ref{versal} and
7.\ref{moduli}. In the Appendix we show how the universal covers of Section  \ref{fund} are organized into families. Theorem 8.\ref{properdiscont} stresses again the importance of the logarithmic type.


\section{Preliminary facts on surfaces with global spherical shells}\label{facts}
In this section we recall some facts on surfaces with global spherical shells.
We refer to \cite{Kato78} and \cite{Dl84} for more details. 

By a {\it surface}
we always mean a compact connected complex manifold of dimension $2$.
A {\it global spherical shell (GSS)} in a surface $S$ 
is the image $\Sigma$ of $S^3$ 
through a  holomorphic imbedding
of an open neighbourhood of $S^3 \subset \mathbb C^2 \setminus \{0\}$
such that $S \setminus \Sigma$ is connected.\\

A minimal surface $S$ with $b_2:=b_2(S)>0$
admits a GSS if it can be obtained by the following method:
Blow up the origin of the unit ball $B$ in $\mathbb C^2$ and choose
a point $p_1$ on the exceptional curve $C_1$. Continue by blowing
up this point and by choosing a further point $p_2$ on the new exceptional
$(-1)$-curve $C_2$. Repeat this process $b_2$ times, i.e. until you
reach the curve $C_{b_2}$ and choose a point $p_{b_2}$ on this curve.
Call $\hat B$ the manifold thus obtained and $\pi:\hat B \to B$ the 
blowing down map. Choose now a biholomorphic map $\sigma:\bar B \to
\sigma(\bar B)$ onto a compact neighbourhood of $p_{b_2}$.
Use finally $\sigma \circ \pi$ to glue together the two components
of the boundary of $\hat B \setminus \sigma(B)$. It is not difficult to see
that in this way one gets a minimal surface $S$ 
with $\pi_1(S)\simeq \mathbb Z$
and $b_2(S)=b_2$. There are exactly $b_2$ rational curves on $S$ which
 form a divisor the dual graph of which 
has one of the three types described in the introduction. 
Besides these rational curves at most one further curve might appear on
 some Enoki surfaces and this curve will be elliptic. One can see a GSS
 as the image of $\sigma(S^3)$ in $S$ via the
 above identification. The universal covering $\tilde S$ of $S$
 is obtained by glueing an infinite number, indexed by $\mathbb Z$,
 of copies of $\hat B \setminus \sigma(B)$ through maps analogous to
 $\sigma \circ \pi$. \\
 
 Conversely, let $S$ be a minimal surface with a fixed GSS $\Sigma$. 
One can fill
 the pseudoconcave end of $S \setminus \Sigma$ holomorphically 
 by $\bar B$ and one obtains a holomorphically convex complex manifold $M$.
 The maximal compact complex analytic set in this manifold is exceptional
 and the blow-down of $M$ is isomorphic to $B$. One shows that this is the
 inverse of the above construction. \\

 In \cite{Dl84}, Dloussky remarked an important object associated
 to this construction, which is the germ $\varphi$ of the holomorphic
 mapping $\pi \circ \sigma:B \to B$ around the origin of $B$.
It is a germ of a contracting mapping and it is shown to completely determine
 the surface $S$ up to isomorphism.
 Conjugating $\varphi$ by an automorphisms of $(B,0)$ does not
 change the isomorphy class of $S$, but to the same surface $S$
 several conjugacy classes of germs may correspond. 
In fact the conjugacy class
 of the contracting germ depends on 
the homotopy class of the GSS $\Sigma \subset S$.\\
 
 In \cite{Fav00}, Favre gave polynomial normal forms for the contracting
 germs and classified them up to biholomorphic conjugacy. In order to obtain
 logarithmic moduli spaces, we shall choose one of these normal forms
 and describe it more closely 
as well as its relation to the maximal reduced divisor
 $D$ of rational curves on the surface. 
A choice of a normal form
 will be suggested by an examination of the fundamental group 
of $S \setminus D$. (Since
$D$ is contractible and does not contain $(-1)$-curves, $S$ is completely 
determined by $S\setminus D$.) 


\section{The fundamental group of $S \setminus D$}\label{fund}

 For a surface $S$ admitting a GSS, one finds using the exponential sequences
 $$ Pic^0(S) 
\simeq H^1(S,\mathbb C^*)\simeq \mathbb C^*\simeq Hom(\pi_1(S),\mathbb C^*).
 $$
 Therefore for each $\lambda \in \mathbb C^*$ 
there is a unique associated flat line bundle 
which we denote by $L_{\lambda}$.
Note that 
the above isomorphisms depend on the choice of a positive generator 
of $\pi_1 (S) \simeq \mathbb Z$. 
We recall also that for the intermediate surface $S$ 
there exists a positive integer $m$, a flat line bundle $L$ 
and an effective divisor $D_m$ 
such that $(K_S \otimes L)^{\otimes m} =\mathcal O_S(-D_m)$,
 cf. \cite{DOT03}, Lemma 1.1. 
The smallest possible $m$ in the above formula is called the 
 {\it index} of the surface $S$ and is denoted by $m(S)$.
By Proposition 1.3 of \cite{DOT03}, 
for each intermediate surface $S$ of index $m$ 
there is a unique intermediate surface $S'$ 
together with a proper map $S' \to S$ 
which is generically finite of degree $m$ such that $S'$ is of index $1$. 
Moreover, on the complements of the maximal reduced divisors $D'$ and $D$ 
the induced mapping $S' \setminus D'\to S \setminus D$ 
is a cyclic unramified covering of degree $m$. \\


 \subsection{Surfaces of index $1$}

In this subsection 
we assume that $S$ is an intermediate surface of index $1$.
Later on we will show how to recover normal forms of contracting maps 
for surfaces of higher index. 
In \cite{DOT01} and  \cite{DOT03} normal forms for the generators of the 
fundamental group of $S\setminus D$ were given. 
We will now recall and further
develop these computations. 

It is shown loc. cit. that the universal cover 
  $\widetilde {S \setminus D}$ of  $S \setminus D$ is isomorphic to 
 $\mathbb C\times  \mathbb H_{l} $, 
where $\mathbb H_l:= \{w \in \mathbb C \mid\Re e(w) <0 \}$ 
is the left half plane. We take $(z,w)$ to be
a system of holomorphic coordinates on  $\mathbb C\times  \mathbb H_{l} $.
The fundamental group of $S \setminus D$ is isomorphic to
 $\mathbb Z \ltimes \mathbb Z[1/k]$ and is generated by the following two 
automorphisms of $\mathbb C\times  \mathbb H_{l} $:
$$\left\{
\begin{array}{lcl}
g_{\gamma}(z,w)&=&(z,w+ 2 \pi i)\\
g(z,w)&=&(\lambda z+a_0+H(e^{-w}),kw),
\end{array}
\right. \leqno{(FG1)}
$$
where $\lambda\in \mathbb C^*$, $H=H(\zeta)= \sum_{m=1}^s a_{m} \zeta^m$
 is a complex polynomial of positive degree 
and $k:=1+\sqrt{|\det M(S)|}\in \mathbb Z, k \geq 2$.
Here $M(S)$ 
denotes the intersection matrix of the $b_2(S)$ rational curves of $S$.
Moreover the complex parameters $\lambda, a_i$ 
are further subject to the following conditions:
$(\lambda-1) a_0=0$ and  $a_{m}= 0$
for all $m>0$ with $k \vert m$, see Theorem 3.1 in \cite{DOT01} and 
Theorems 4.1 and 4.2 in \cite{DOT03}.
The case $\lambda=1$ occurs precisely when 
$S$ admits a nontrivial holomorphic vector field. Note also that
the generator $g_{\gamma}$ corresponds to some loop $\gamma$ around $D$,
whereas $g$ induces a generator of $\pi_1(S)$. It is claimed in \cite{DOT01}
and \cite{DOT03} that this normal form is unique if one preserves the 
generator $g_{\gamma}$. We shall see later that this claim is valid only
up to the
action of some finite group on the space of parameters 
given by the coefficients of $H$.\\

 We first show how the normal form changes if one allows
a different choice of the loop $\gamma$. 

We start by replacing the normal form (FG1) as in \cite{DOT01}, p. 659 and \cite{DOT03}, p. 307,
in order to obtain a new normal form which will lead to the contracting germ:
$$\left\{
\begin{array}{lcl}
g_{\gamma}(z,w)&=&(z,w+ 2 \pi i)\\
g(z,w)&=&(\lambda z+a_0+Q(e^{-w}),kw),
\end{array}
\right. \leqno{(FG2)}
$$
where $Q=Q(\zeta):= \sum_{m=l}^{\sigma} b_m \zeta^m$ is a new polynomial
such that $b_{\sigma}=1$, $k \nmid \sigma$, $l := [\sigma/k] +1$ and 
$\gcd \{k,m \mid b_m \neq 0\} =1$. As before $(\lambda-1)a_0=0$. 
The normalisation $b_{\sigma}=1$ is obtained by a conjugation
$ (z,w) \mapsto (\alpha z , w)$ with a suitable $\alpha$.\\

We introduce now a numerical invariant of (FG2) which will be shown to be 
equivalent to the dual graph of $D$, see Section \ref{blow-up}.

\begin{Def} \label{type}
Let $k$ be fixed and let $Q$ be a polynomial as in (FG2).
Define inductively the
following finite sequences of integers $\sigma=:n_1>....>n_t \geq l$ 
and $k> j_1> j_2...>j_t=1$  by:
\begin{itemize}
\item[i)] $n_1:=\sigma$, $j_1:=\gcd(k,n_1)$;
\item[ii)] $n_\a:=\max\{i<n_{\a -1} \mid b_i \neq 0, \gcd(j_{\a -1},i) < j_{\a -1}\},\\
j_\a: = \gcd(k,n_1,...,n_{\a })=\gcd(j_{\a -1}, n_\a)  $;
\item[iii)] $1=j_t:=\gcd(k,n_1,...,n_{t -1},n_t) <
\gcd(k,n_1,...,n_{t -1})$.
\end{itemize}
We call $(n_1,...,n_t)$ the {\rm type} of (FG2) and $t$ the {\rm length of the type}.
If $t=1$ we say that (FG2) is of {\rm simple type}.
\end{Def}

We remark now that for any divisor $d >1$ of $k$ the following are also generators
of $\pi_1(S \setminus D) \subset {\rm Aut}_{\mathcal O}(\mathbb C \times \mathbb H_l)$:
$$\left\{
\begin{array}{lcl}
g_{\gamma}^d (z,w)&=&(z,w+ 2  \pi i d)\\
g(z,w)&=&(\lambda z+a_0+Q(e^{-w}),kw),
\end{array}
\right. \leqno{(FG3)}
$$  
since $(g^{-1} \circ g_{\gamma}^d \circ g)^{\frac{k}{d}} =g_{\gamma}$.\\
In order to reestablish (FG2) we perform the following 
conjugations on (FG3).
First conjugate with $(z,w) \mapsto (z,dw)$. The generators become
$$\left\{
\begin{array}{lcl}
  (z,w)& \mapsto &(z,w+ 2  \pi i )\\
 (z,w)&\mapsto &(\lambda z+a_0+Q(e^{-dw}),kw),
\end{array} 
\right.  \leqno{(1)}
$$  
This is not in normal form (FG2); in fact there are two cases.
If $k \nmid d \sigma$ one has to go back by the inverse conjugation procedure
to the initial (FG2). 
At the end of this process the loop $\gamma $ remains unchanged.\\
If $k \mid d \sigma$ 
take $p:= \max\{i \in \mathbb N \mid i >1,\  k\! \mid \! j_{i-1}d \}$.
Now let  
$$\Phi_1(z,w) = (z+ \sum_{m=n_p +1}^{\sigma} b_m e^{-\frac{mdw}{k}},w) $$ and apply the conjugation $\Phi_1 \circ (1) \circ \Phi_1^{-1}$.
We obtain the new generators
$$\left\{
\begin{array}{lcl}
(z,w)  &   \mapsto &  (z,w+ 2  \pi i )\\
(z,w) &  \mapsto &  (\lambda z+a_0+b_{n_p} {\tilde Q_{d}}(e^{-w}), kw)
\end{array} 
\right.  \leqno{(2)}
$$  
where
$$\tilde Q_d (\zeta)
:= b_{n_p}^{-1} \Big( Q(\zeta^d )+\lambda \sum_{m=n_p +1}^{\sigma}
 b_m \zeta^{\frac{md}{k}} -\sum_{m=n_p +1}^{\sigma} b_m \zeta^{md} \Big).$$
 
 The conjugation $\Phi_2 \circ (2) \circ \Phi_2^{-1}$ by 
$\Phi_2(z,w) = (b_{n_p}z,w)$ leads us to
$$\left\{
\begin{array}{lcl}
(z,w)  &   \mapsto &  (z,w+ 2  \pi i )\\
(z,w) &  \mapsto &  (\lambda z+b_{n_p}^{-1}a_0+ {\tilde Q_{d}}(e^{-w}), kw).
\end{array} 
\right.  \leqno{(3)}
$$  
The idea is to bring (3) in normal form (FG2) presumably with the new type
$(\tilde n_1,...,\tilde n_t)
=(dn_p,...,dn_t,\frac{dn_1}{k},...,\frac{dn_{p-1}}{k}).$
We have $\gcd(dn_p,k)=\\ \gcd(dn_p,dj_{p-1},k)$ because $\gcd(dj_{p-1},k)=k$ and
 further $\tilde j_1=\gcd(dn_p,k)=\gcd(dn_p,dj_{p-1},k)
=\gcd(dj_{p},k)<k$. In the same way one obtains 
$\gcd(dm,k)\\
=\gcd(dj_{p},k)=   \tilde j_1$ for $n_{p+1}<m<n_p$, 
  $\tilde j_2=  \gcd(dn_{p+1},dn_p,k)=$\\ $\gcd(dn_{p+1},dj_{p},k)
=\gcd(dj_{p+1},k)$ and so on. The last in this series is
$\gcd(dn_{t},dj_{t-1},k)$
$=\gcd(dj_{t},k)= \gcd(d,k)=d$.
The translated exponents $\frac{dn_1}{k},...,\frac{dn_{p-1}}{k}$ now
give  $${\tilde j}_{t-p+2}=\gcd(d,{\tilde n}_{t-p})= \gcd(\frac{dn_1}{k},d,k) =\gcd(\frac{dn_1}{k},d)=\gcd(\frac{dj_1}{k},d)=\frac{dj_1}{k}<d
$$
because $k \nmid n_1$, 
 $${\tilde j}_{t-p+3}=\gcd(\frac{dn_2}{k},\frac{dn_1}{k},k) =\gcd(\frac{dn_2}{k},\frac{dj_1}{k})=\frac{dj_2}{k},
$$
 and so on, down to ${\tilde j}_t=\frac{dj_{p-1}}{k}$.
Set  $d':= \frac{dj_{p-1}}{k}$. If $d'=1$, then $\tilde Q_d$ is in normal form and
we are ready. If not all its "active" exponents $m$, i.e. with $b_m \neq 0$,
are divisible by $d'$. Then
we perform the same conjugations which lead us
from (FG3) to (1) in inverse order and with $d'$ 
instead of $d$. We get 
$$\left\{
\begin{array}{lcl}
(z,w)  &   \mapsto &  (z,w+ 2  \pi i )\\
(z,w) &  \mapsto &  (\lambda z+b_{n_p}^{-1}a_0+ {\tilde Q_{d}}(e^{-\frac{w}{d'}}), kw).
\end{array} 
\right.  \leqno{(4)}
$$  

The type of this normal form is 
$$\frac{1}{d'} (dn_p,...,dn_t,\frac{dn_1}{k},...,\frac{dn_{p-1}}{k}).$$
It is clear that we get exactly $t$ (FG2)-normal forms through this kind of transformation
corresponding to the following choices for $d$:
$\frac{k}{j_1},...,\frac{k}{j_t}$, the last one giving our starting type again. 
In these cases we also have $d'=1$, hence the type transformations are of the form:
$$ (n_1,...,n_t) \mapsto (\frac{kn_{p+1}}{j_p},...,\frac{kn_t}{j_p},\frac{n_1}{j_p},...,
\frac{n_p}{j_p}).$$
This amounts to an action of $\mathbb Z/(t)$ on the set of types of length $t$.\\

Notice moreover that we also get a biholomorphic polynomial map $\tau_d$ between
the parameter spaces for the coefficients of $Q$ in (FG2) and those of
$\tilde Q_d$ in (4).


\subsection{Surfaces of higher index}\label{highindex}

We consider now the normal form (FG2) for a surface $S$ of index $1$
and a positive integer $q$.
We say that the set of parameters $a_0,b_l,...,b_{\sigma}$ of the normal form
has the property $(I_q)$ if
$$a_0=0$$ 
and
$$q | \gcd \{k-1; \ m-m', \ b_mb_{m'}\neq 0\}.$$
Suppose now that this is the case and conjugate $g_{\gamma}$ and $g$ by
$$\varphi (z,w):=(e^{\frac{\sigma 2\pi i}{q}}z,w-\frac{2\pi i}{q}).$$
We get
$$\varphi\circ g_{\gamma}\circ \varphi^{-1}=g_{\gamma},$$
$$\varphi\circ g\circ \varphi^{-1}=g_{\gamma}^{\frac{k-1}{q}}\circ g,$$
showing that $\varphi$ lies in the normalizer of the fundamental group $\pi_1(S\setminus D)$.
This leads to an action of $\Z/(q)$ on $S\setminus D$ and to a quotient
$S'\setminus D'$, where $S'$ is a surface of possibly higher index. 
We will see
later that every surface of index $q$ arises in this way, see Remark 4.\ref{automorph}. 
Looking now at a divisor $d$ of $k$ we notice that the transformation
$Q \stackrel{\tau_d}{\mapsto} \tilde Q_d$ of the preceding subsection preserves the property $(I_q)$. This happens because
$d$ and $q$ are relatively prime. Our next step will be to translate these properties for 
the associated contracting germs of holomorphic maps.


\section{Contracting germs of holomorphic mappings associated to intermediate surfaces}
\label{germs}

In \cite{Fav00} Favre gave normal forms for contracting rigid germs of holomorphic maps
$(\C^2,0)\to(\C^2,0)$ and completely characterized those germs which give rise to surfaces with global spherical shells. We retain the following result which will be essential for our study.

\begin{Th}[Favre \cite{Fav00}, see also \cite{DOT03}]\label{favre} 
\noindent Every intermediate surface is associated to a polynomial germ in the following normal form:
$$  \varphi(z, \z) := (\l \z^s z +P(\z)+c\z^{\frac{sk}{k-1}}, \z^k),   \leqno{(CG)}$$
where $k,s\in\Z$, $k>1$, $s>0$, $\l\in\C^*$, $P(\z):= c_j \z^j+c_{j+1}\z^{j+1}...+c_s\z^s$ is a
complex polynomial, $0<j<k$, $j \leq s$, $c_j=1$, $c_{\frac{sk}{k-1}}:=c\in \C$ with $c=0$ whenever
${\frac{sk}{k-1}}\notin \Z$ or $\l\neq 1$ and $\gcd \{k,m \mid c_m \neq 0\} =1$.\\
Moreover, two polynomial germs in normal form (CG) $\varphi$ and 
$$\tilde{\varphi}:= (\tilde{\l} \z^{\tilde{s}} z +\tilde{P}(\z)+\tilde c \z^{\frac{\tilde s \tilde k}{\tilde k-1}}, \z^{\tilde{k}})$$ are conjugated if and only if there exists $\e\in\C$ with $\e^{k-1}=1$ and
$\tilde k=k$, $\tilde s=s$,
$\tilde{\l}=\e^s \l$, $\tilde P(\z)=\e^{-j}P(\e \z)$, $\tilde c=\e^{\frac{sk}{k-1}-j}c$.
\end{Th}

\begin{Rem}
 Intermediate surfaces of index one correspond precisely to germs $\varphi$ in normal form (CG) such that
$(k-1)\vert s$.
\end{Rem}
In order to see this we look at a surface $S$ of index one and at the normal form (FG2) 
of the fundamental group of $S\setminus D$. We conjugate (FG2) by 
$(z,w)\mapsto (e^{lw}z,w)$ and get the following form for the generators of
$\pi_1(S\setminus D)$:
$$\left\{
\begin{array}{lcl}
(z,w)  &   \mapsto &  (z,w+ 2  \pi i )\\
(z,w) &  \mapsto &  (\lambda e^{l(k-1)w}z+e^{lkw}a_0+P(e^{w}), kw)
\end{array} 
\right.  
$$  
where $P$ is the polynomial defined by $P(\z)=\z^{lk}Q(\z^{-1})$.
Thus $P(\z)=\sum_{m=lk-\s}^{l(k-1)} c_m \zeta^m$, whith $c_m=b_{lk-m}$. As shown in 
\cite{DOT01} and  \cite{DOT03} the surface $S$ is then associated to the polynomial germ
$$(z,\zeta)\mapsto (\l \z^s z +P(\z)+c\z^{\frac{sk}{k-1}}, \z^k)$$
which is in normal form (CG). Here we have set $s=l(k-1)$ and $c=a_0$.

\begin{Rem}\label{automorph}
 Every intermediate surface of index higher than  one may be constructed as 
in section \ref{highindex}.
\end{Rem}
Indeed, let $\varphi '$ be a polynomial germ in normal form (CG) associated to
a surface $S'$ of higher index,
$$\varphi '(z,\z)=(\l\z^{s'}z+\sum_{m=j'}^{ s'} c'_m \zeta^m,\z^k),$$
with $(k-1)\nmid s'$ and let $q$ be some positive divisor of $k-1$ such that
$(k-1)$ divides $qs'$. Set $r:=\lfloor \frac{qj'}{k} \rfloor$, 
$s:=qs'-r(k-1)$, $j:=qj'-rk$, $P(\z):=\sum_{m=j'}^{ s'} c'_{m} \zeta^{qm-rk}$ and
$$\varphi(z,\z):= (\l\z^{s}+P(\z),\z^k).$$ 
One checks now easily that the polynomial germ $\varphi$ is in normal form (CG),
corresponds to a surface $S$ of index one and admits an automorphism of the form
$$(z,\z)\mapsto (\e^{-r}z,\e\z),$$
where $\e$ is a primitive root of unity of order $q$. 
This automorphism lifts to a conjugation of a corresponding normal
form (FG2) for $S$, which has property $(I_q)$ as in section \ref{highindex}.
In fact the covering map $S\setminus D\to S'\setminus D'$ is induced by
$(z,\z)\mapsto (\zeta^r z,\z^q)$. \\
One also sees that the index of $S'$ is the least possible $q$ allowed in the above construction;
more precisely, one gets

\begin{Rem}
The index of an intermediate surface associated to a polynomial germ
$\varphi(z, \z)= (\l \z^s z +P(\z)+c\z^{\frac{sk}{k-1}}, \z^k)$ in normal form (CG) equals
$$\frac{k-1}{\gcd (k-1,s)}.$$
\end{Rem} 

\begin{Rem}(Cf. \cite{DOT01}.)\\
 The intermediate surfaces admitting non-trivial holomorphic vector fields  are precisely those 
associated to polynomial germs  $$  \varphi(z, \z) := (\l \z^s z +P(\z)+c\z^{\frac{sk}{k-1}}, \z^k)$$
in normal form (CG) with $(k-1)\vert s$ and $\l=1$.
\end{Rem}

We now define the type of a polynomial germ in normal form (CG) in a similar way we did it for (FG2).

\begin{Def} For fixed $k$ and $s$ and for a polynomial germ 
$  \varphi(z, \z) := (\l \z^s z +P(\z)+c\z^{\frac{sk}{k-1}}, \z^k)$ in normal form
(CG) with $P(\z):= \z^j+c_{j+1}\z^{j+1}+...+c_s\z^s$ we define inductively the
following finite sequences of integers $j=:m_1<....<m_t \leq s$ 
and $k> i_1> j_2...>i_t=1$  by:
\begin{itemize}
\item[i)] $m_1:=j$, $i_1:=\gcd(k,m_1)$;
\item[ii)] $m_\a:=\min\{m>m_{\a -1} \mid c_m\neq 0, \gcd(i_{\a -1},m) < i_{\a -1}\},
i_\a: = \gcd(k,m_1,...,m_{\a })=\gcd(i_{\a -1}, m_\a)  $;
\item[iii)] $1=i_t:=\gcd(k,m_1,...,m_{t -1},m_t) <
\gcd(k,m_1,...,m_{t -1})$.
\end{itemize}
We call $(m_1,...,m_t)$ the {\rm type} of $\varphi$ and $t$ the {\rm length of the type}.
If $t=1$ we say that $\varphi$ is of {\rm simple type}.
\end{Def}

The length of the type of $\varphi$ is $1$, i.e. $\varphi$ is of simple type, if
and only if $k$ and $j$ are relatively prime, $\gcd (k,j)=1$. 

We also set
$\varepsilon(k,m_1,...,m_t,s):=\lfloor \frac{m_2-m_1}{i_1} \rfloor+
\lfloor \frac{m_3-m_2}{i_2} \rfloor+...+
\lfloor \frac{m_t-m_{t-1}}{i_{t-1}} \rfloor+
s-m_t$. It is the number of coefficients of $P$ whose vanishing or non-vanishing does not affect the type.

The type is obviously preserved by the conjugations appearing in Theorem 4.\ref{favre}.

For fixed $k$, $s$ and fixed type $(m_1,...,m_t)$ we consider the following obvious
 parameter spaces for the coefficients $(\l,c_{j+1},...,c_s,c)$ appearing in 
(CG):
\begin{itemize}
\item when $(k-1)$ does not divide $s$ we take 
$$U_{k,s,m_1,...,m_t}=\C^*\times(\C^*)^{t-1}\times\C^{\varepsilon(k,m_1,...,m_t,s)},$$ 
\item in case $(k-1)\mid s$
$$U_{k,s,m_1,...,m_t}^{\l\neq1,c=0}=\C\setminus \{0,1\}\times(\C^*)^{t-1}\times\C^{\varepsilon(k,m_1,...,m_t,s)},$$
 $$U_{k,s,m_1,...,m_t}^{\l=1}=(\C^*)^{t-1}\times\C^{\varepsilon(k,m_1,...,m_t,s)}\times\C.$$
\end{itemize}
In the second case we have considered separate parameter 
spaces for surfaces without,
respectively with, holomorphic vector fields. We will come back later to this point.

By the discussion of the previous section we see that for each $p\in\{1,...,t\}$ choosing 
$d=\frac{k}{i_p}$ gives us a transformation on types and a biholomorphic map $\tau_d$ between
corresponding parameter spaces. Surfaces associated to germs of one type correspond via
 $\tau_d$ to isomorphic surfaces associated to germs of the
transformed type.


\section{Decomposition of germs}

\begin{Def} A polynomial germ
$$  \varphi(z, \z) := (\l \z^s z +P(\z)+c\z^{\frac{sk}{k-1}}, \z^k)$$
in normal form (CG) is said to be 
in {\rm pure normal form} if $c=0$.
We say that a polynomial germ $\varphi:(\C^2,0)\to(\C^2,0)$ is in {\rm modified normal form} 
if
$$  \varphi(z, \z) := (\l \z^s z +P(\z)+c\z^{kn}, \z^k)  $$
is such that the germ
$$  \sideset{_{p}}{}{\f}(z, \z):= (\l \z^s z +P(\z), \z^k)  $$
is in pure normal form (CG), $n\in \Z$,  $\frac{s}{k}< n \leq \frac{s}{k-1}$ and $c\in \C$. 
In this case the {\rm type} of $\varphi$ is by definition the same as the type of
the {\rm  purified germ} $\sideset{_p}{}{\f}$.
\end{Def}

\begin{Rem}\label{modif}
 The pure normal form is a special case of the normal form (CG). The  normal form (CG) is a special case of the 
modified normal form.
A germ  $  \varphi(z, \z) := (\l \z^s z +P(\z)+c\z^{kn}, \z^k)$  in modified normal form is  
in normal form (CG) if and only if either $n=\frac{s}{k-1}$ and $\l=1$ or $c=0$. 
When this is not the case
$\varphi$ is conjugated to some germ $\tilde{\varphi}:= (\l \z^{s} z +\tilde{P}(\z), \z^{k})$ 
in pure normal form and of the same type as $ \varphi$.
\end{Rem} 
For suppose $c\neq 0$. We have to consider two cases. 

When $ \frac{s}{k}<n<\frac{s}{k-1}$ we conjugate 
$\varphi$ by
$  (z, \z)\mapsto (z - \l^{-1}c\z^{kn-s}, \z)  $ to  get the new germ 
$$  (z, \z)\mapsto (\l \z^s z +P(\z) +\l^{-1}c\z^{k(kn-s)}, \z^k).  $$
By assumption we have $kn-s<n$. If $k(kn-s)\leq s$ we are ready since 
$  (z, \z)\mapsto (\l \z^s z +P(\z) +\l^{-1}c\z^{k(kn-s)}, \z^k)  $ 
is in pure normal form. If not we can work with $kn-n$ instead of $n$ and continue
to conjugate until the exponent of the supplementary $\z$-term doesn't exceed $s$ any more.
Remark that this exponent remains a multiple of $k$ and thus leaves the type of $\varphi$ intact.

When $ n=\frac{s}{k-1}$ and $\l\neq 1$ a conjugation by 
$  (z, \z)\mapsto (z - \frac{c}{\l-1}\z^{n}, \z)  $ leads $\varphi$ directly to the pure normal form
$$  (z, \z)\mapsto (\l \z^s z +P(\z), \z^k).  $$

\begin{Prop}
 Let $  \varphi(z, \z) := (\l \z^s z +P(\z)+c\z^{kn}, \z^k)$ be a germ 
 in modified normal form and whose type has length
$t$. If $\varphi$ is in pure normal form,
then $\varphi$ admits a canonical decomposition $\varphi=\varphi_1\circ...\circ\varphi_t$ into $t$ polynomial germs
in pure normal form and of simple type.  If  $ \varphi$ is not in pure normal form,
then $\varphi$ admits a canonical decomposition $\varphi=\varphi_1\circ...\circ\varphi_t$ into $t$ polynomial germs in modified
normal form and of simple type out of which $t-1$ are pure.
Moreover $\sideset{_p}{}{\f}=\sideset{_p}{_1}{\f}\circ...\circ\sideset{_p}{_t}{\f}$.
The canonical choice of the decomposition $\varphi=\varphi_1\circ...\circ\varphi_t$ and the types 
of the factors $\varphi_1$,...,$\varphi_t$ are determined by the type of 
$\varphi$. 
\end{Prop} 
\begin{proof}
   We first deal with the pure case.

Let $$ \varphi_1(z, \z):= (\l_1 \z^{s_1} z +P_{1}(\z), \z^{k_1})$$ 
and
$$ \varphi_2(z, \z):= (\l_2 \z^{s_2} z +P_{2}(\z), \z^{k_2})$$
be two germs in pure normal form with invariants $(k_1,j_1,s_1)$ resp. $(k_2,j_2,s_2)$ of types
$$(j_1=m_1^{(1)},m_2^{(1)},...,m_{t_1}^{(1)})$$
and
$$(j_2=m_1^{(2)},m_2^{(2)},...,m_{t_2}^{(2)})$$ respectively.
Their composition 
$$\varphi_1 \circ \varphi_2 (z,\z) =
( \l_1  \l_2 \z^{s_1 k_2 +s_2} z + P_1(\z^{k_2})+ \l_1 \z^{s_1 k_2} P_2(\z),\z^{k_1 k_2})$$

is again in pure normal form with invariants $(k=k_1 k_2, j= j_1 k_2,s= s_1 k_2 +s_2)$
and has type
$$(m_1^{(1)}k_2,m_2^{(1)} k_2,...,m_{t_1}^{(1)} k_2, s_1 k_2+m_1^{(2)},  s_1 k_2+m_2^{(2)},..., s_1 k_2+m_{t_2}^{(2)})$$
with length $t=t_1 +t_2$.

Conversely, let 
$$\varphi(z, \z) :=  (\l \z^s z +P(\z), \z^k),$$
where $$P(\z):=\z^j+c_{j+1}\z^{j+1}...+c_s\z^s$$
be a germ in pure  normal form and suppose that 
$$d:=\gcd(j,k)>1.$$
Define $k_1:=k/d,\ j_1:=j/d$ and $$s_1:= \max\{m/d \mid  
\tilde m /d \in \mathbb Z \ {\rm or} \ c_{\tilde m} =0 \ \ {\rm for \ all} \ \  \tilde m =j,...,m  \}.$$
It is clear that $0< j_1< k_1$ and $j_1 \leq s_1$ . Furthermore let $k_2:=d$, 
$$j_2:= - s_1 k_2 + \min\{m \mid m > s_1 k_2 \ \ {\rm and} \ \ c_m \neq 0\}$$
and $s_2:=s-s_1 k_2$. It is clear that $0< j_2 < k_2  $ and $j_2 \leq s_2$.
We also put $\l_1:= c_{s_1 k_2+j_2} \neq 0$ and $\l_2:= \l/\l_1$. Finally
let 
$$\varphi_1(z,\z) := (\l_1 \z^{s_1} z + \sum_{m=j}^{s_1 k_2} c_m \z^{m/d}, \z^{k_1}),$$
$$\varphi_2(z,\z) := (\l_2 \z^{s_2} z +    \sum_{m=s_1 k_2+j_2}^{s}\l_{1}^{-1} c_m \z^{m-s_1 k_2}, \z^{k_2}).$$
One verifies directly that $\varphi_1 \circ \varphi_2 = \varphi$, that $\varphi_1 $, $\varphi_2 $ are in pure normal form and that
$\varphi_1 $ is of simple type. Note also that the types of $\varphi_1 $ and $\varphi_2 $ are determined by the type of
 $\varphi$. 
Repeating this procedure if necessary with $\varphi_2$ one gets a decomposition
of $\varphi$ into germs in pure normal form and of simple type in a canonical way.

Take now $  \varphi(z, \z) = (\l \z^s z +P(\z)+c\z^{kn}, \z^k)  $
a germ in modified normal form, $  \sideset{_{p}}{}{\f}(z, \z)= (\l \z^s z +P(\z), \z^k)  $
and $ \varphi_1(z, \z):= (\l_1 \z^{s_1} z +P_{1}(\z), \z^{k_1})$,
$ \varphi_2(z, \z):=  (\l_2 \z^{s_2} z +P_{2}(\z), \z^{k_2})$
in pure normal form such that $  \sideset{_{p}}{}{\f}=\varphi_1\circ\varphi_2$.
We shall modify $\varphi_1$ or $\varphi_2$ in order to get $\varphi$ by composition.

If $n\leq \frac{s_1}{k_1-1}$, 
then $  \frac{s_1}{k_1}< n$ and 
the germ 
$$\tilde \varphi_1(z, \z):= (\l_1 \z^{s_1} z +P_{1}(\z) +c\z^{k_1n}, \z^{k_1})$$
is in modified normal form and $\tilde \varphi_1 \circ \varphi_2 = \varphi$.
If not, take $n_2:=k_1n-s_1$ and 
$$\tilde \varphi_2(z, \z):= (\l_2 \z^{s_2} z +P_{2}(\z) +c\l_1\z^{k_2n_2}, \z^{k_2}).$$
Then $\varphi_1 \circ \tilde \varphi_2 = \varphi$. The inequalities
$\frac{s_1}{k_1-1}<n\leq\frac{s}{k-1}$ imply $n_2< \frac{s_2}{k_2-1}$ 
and thus $\tilde \varphi_2 $ will be in modified normal form.

Note that in case $n=\frac{s_1}{k_1-1}=\frac{s}{k-1}=\frac{s_2}{k_2-1}$
both decompositions $ \varphi=\tilde \varphi_1 \circ \varphi_2 =\varphi_1 \circ \tilde \varphi_2 $
are possible but we chose the first $\tilde \varphi_1 \circ \varphi_2$ as the canonical one.
\end{proof}


\section{The blow-up sequence and the Dloussky sequence} \label{blow-up}

In this section we calculate for a given germ in modified 
normal form the configuration of the rational curves
on the associated intermediate surface. This is equivalent
to giving the {\it Dloussky sequence} of the surface, see \cite{Dl84}, pp. 37,
which in turn will be computed by the sequence of blow-ups
for the given germ. We recall that in the Dloussky sequence
each entry represents a rational curve with the negative of its self-intersection  
and the order in the sequence is given by the order of creation of the curves in the blow-up
process. A Dloussky sequence for an intermediate surface
is called {\it simple} if it is of the form

\begin{eqnarray*}
[{\rm DlS}]&=& [\a_1 +2,\underbrace{2,..,2}_{\a_1 -1},...,\a_q +2,\underbrace{2,..,2}_{\a_q -1},\a_{q+1} +1,\underbrace{2,..,2}_{\a_{q+1} -2},\underbrace{2,..,2}_{m}]=\\
&=&[s_{\a_1},...,s_{\a_{q}},s_{\a_{q+1}-1},r_{m}], 
\end{eqnarray*}

with $q \geq 0$, $\a_i \geq 1$ for $ 1 \leq i \leq q$, $\a_{q+1} \geq 2$ and $m \geq 1$. So the shortest possible sequence here appears for
$q=0$ and $m=1$ and is of the form $[3,2]$. A general Dloussky sequence is of the form
$[{\rm DlS}_1,....,{\rm DlS}_N]$, where $[{\rm DlS}_j], j=1,...,N$ are simple Dloussky sequences. The dual graph of the divisor $D$
on an intermediate surface with known Dloussky sequence
is now constructed in the following way: The entries of the 
sequence represent the knots of the graph. An entry with value $\a$ 
is connected with the entry following $\a -1$ places after it at
the right hand (with the entries in cyclic order!). 

\medskip

\begin{Ex} a)
The Dloussky sequence $[3,4,2,2]$ produces the graph

\medskip

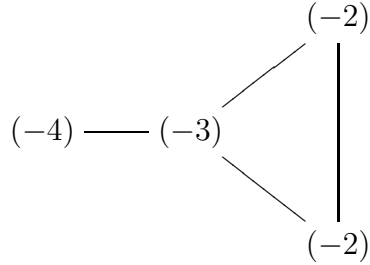
\xymatrix{ 
&&&&  & & {(-2)} \\ 
&&&&  {(-4)}\ar@{-}[r]   &  {(-3)} \ar@{-}[ur] \ar@{-}[dr] &\\
&&&&  & &   {(-2)} \ar@{-}[uu] } 
\medskip

b) The Dloussky sequence $[3,2,4,2,2,2]$ produces the graph

\medskip
 
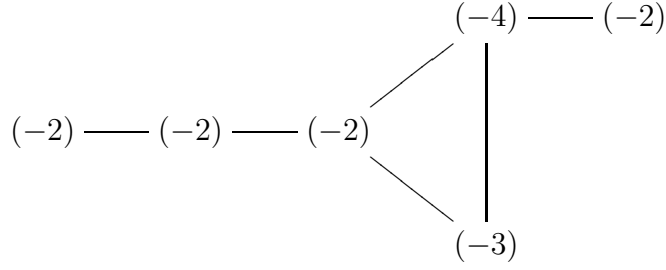
\xymatrix{ 
&&&&  & & {(-4)} \ar@{-}[r]& (-2) \\ 
&&&{(-2)}\ar@{-}[r] &  {(-2)}\ar@{-}[r]   &  {(-2)} \ar@{-}[ur] \ar@{-}[dr] &&\\
&&&&  & &   {(-3)} \ar@{-}[uu] &}
\medskip
\end{Ex}

Note that the Dloussky sequence induces an orientation on the cycle of the dual graph. We will call such a graph a {\em directed dual graph}. Note also that the dual graph of the cycle of rational curves on a surface with a GSS has a natural orientation given by the order of creation of the curves in the blowing-up process. After glueing as described in Section \ref{facts} this orientation may be recovered by looking at the pseudoconvex side of a GSS in the compact surface $S$.

We start with the description of the blow-up sequence in the case of a germ $\varphi $ of simple type. This is the sequence of blow-ups which leads to the maps  $\pi:\hat B \to B$,
$\sigma:\bar B \to
\sigma(\bar B)$ such that $\varphi=\pi\circ\sigma$ as in Section \ref{facts}.
 
Let 
 $$\varphi(z, \z) := (\l \z^s z + \z^j+ 
c_{j+1}\z^{j+1}...+c_s\z^s +c\z^{kn} ,\z^k)=$$
$$= (\l \z^s z +P(\z)+c\z^{kn}, \z^k)= (\z^j A(z,\z), \z^k)= (\z^j A,\z^k),$$
where $$P(\z):=\z^j+c_{j+1}\z^{j+1}...+c_s\z^s$$ and
$$A(z,\z):= 
\l \z^{s-j} z + 1+c_{j+1}\z +...+c_s\z^{s-j}+c\z^{kn-j},\ A(0,0)=1.$$
be a germ in modified normal form and suppose that 
$d:=\gcd(j,k)=1$, that is, $\varphi$ is of simple type.
Recall that $0<j<k$, $j \leq s$ and let $r=s-j \geq 0$.
We shall calculate the blow-up sequence of $\varphi$.
For this we shall need the division algorithm for $k$ and $j$
which we note as follows:

\begin{eqnarray*}
 k &  = & {\a_1 j  +  \b_1} \\
 j & = & {\a_2 \b_1  +  \b_2} \\
 \b_1 & = & {\a_3 \b_2  +  \b_3} \\
&\vdots & \\
 \b_{q-2}  &= &  \a_q \b_{q-1}  +  1 \\
 \b_{q-1} & = &  \a_{q+1}\cdot 1   +  0   
\end{eqnarray*}

with the convention $k=\b_{-1}, j=\b_0$.\\
We use the two standard blow up coordinate charts 
$\eta:(u,v) \to (uv,v)$ and $\eta':(u',v') \to (v',u'v')$ for this procedure.
Remark that in these charts the exceptional curve is given by 
$v=0$ respectively by $v'=0$.

$$
 \begin{xy}
  \xymatrix{
      (\zeta^j A, \zeta^k)             
      &  (\zeta^{k-j}A^{-1},\zeta^jA)  \ar[l]_{   \eta'}  \\
             & \\
(\zeta^{(k-\a_1j)=\b_1}A^{-\a_1},\zeta^jA) \ar[ruu]^{  \eta^{\a_1 -1}}           & (\zeta^{j-\b_1}A^{a_1+1},\zeta^{\b_1}A^{-a_1})  \ar[l]_{   \eta'}  \\
            & \\
(\zeta^{(j-\a_2 \b_1)=\b_2}A^{ \a_1\a_2 +1},\zeta^{\b_1}A^{-\a_1}) \ar[ruu]^{  \eta^{\a_2 -1}}          & (\zeta^{\b_1-\b_2}A^{\cdots},\zeta^{\b_2}A^{\cdots})  \ar[l]_{   \eta'}  \\
\vdots   &  \vdots    \\ 
 (\zeta^{(\b_{q-2}-\a_q \b_{q-1})=\b_q=1}A^{\cdots},\zeta^{\b_{q-1}}A^{\cdots})             & (\zeta^{\b_{q-1}-1}A^{\cdots},\zeta A^{\cdots})  \ar[l]_{   \eta'}  \\
            & \\
(\zeta A^{\cdots},\zeta A^{\cdots}) \ar[ruu]^{  \eta^{\a_{q+1} -2}} &           ( A^{\cdots}-1,\zeta A^{\cdots})  \ar[l]_{  (((u+1)v,v) \leftarrow (u,v)}  \\
     &\\
 (C\lambda z+ ...,\zeta A) =:\sigma(z,\zeta) \ar[ruu]^{\tilde\eta}            & }
\end{xy}
$$
   
In the above sequence of blow-ups all the powers of $A$ indicated
by $\cdots$ are non-zero integers. The map $\sigma$
is a germ of a biholomorphism at the origin of $(\C^2,0)$
with the property that the inverse image $\s^{-1}(C \cap U)$
of the intersection of the last created rational curve $C$ with
an open neighbourhood $U$ of $\s(0,0)$ is given by $\{ \z =0\}$.
Furthermore $\tilde \eta$ is the composition of $s-j$
blowups in the form $(u,v) \mapsto ((u+const.)v,v)$.\\

The construction shows that the Dloussky sequence of the associated surface is simple and given by

\begin{eqnarray*}
[{\rm DlS}]&=& [\a_1 +2,\underbrace{2,..,2}_{\a_1 -1},...,\a_q +2,\underbrace{2,..,2}_{\a_q -1},\a_{q+1} +1,\underbrace{2,..,2}_{\a_{q+1} -2},\underbrace{2,..,2}_{s-j+1}]=\\
&=&[s_{\a_1},...,s_{\a_{q}},s_{\a_{q+1}-1},r_{s-j+1}], 
\end{eqnarray*}
see \cite{Dl84}, pp. 37.

In particular, the second Betti number of the surface
is $ b_2=(-1+\sum_{i=1}^{q+1} \a_i)   +(s-j+1)= (\sum_{i=1}^{q+1} \a_i )  +(s-j)$ which is the {\it length} of $[\rm DlS]$.

\begin{Rem}
The above construction shows that the obtained Dloussky sequence
is independent of the germ being in normal form, pure normal form
or modified normal form. 
\end{Rem}

For the general case let $\varphi =\varphi_1 \circ \ff_2 \circ...
\circ \ff_N $ be the decomposition of a germ $\ff$ into germs of simple type and
${\rm DlS}$ resp. ${\rm DlS}_i$ the associated 
Dloussky sequences of $\ff$ resp. $\ff_i,i=1,...,N$. An easy calculation similar to the  above
one 
shows that  $$[{\rm DlS}]=[{\rm DlS}_1,....,{\rm DlS}_N],$$
i.e. the operations of composition of germs and concatenations
of Dloussky sequences are compatible.

We note in conclusion that the following three objects associated to an intermediate surface are algorithmically computatble from one another: the directed dual graph of the rational curves, the Dloussky sequence, the type of a contracting germ in modified normal form.

\section{Versal families and moduli spaces}\label{secmoduli}

We have seen in section \ref{germs} that for fixed $k$, $s$ and fixed type $(m_1,...,m_t)$
we get parameter spaces for germs in normal form (CG)
$$U_{k,s,m_1,...,m_t}=\C^*\times(\C^*)^{t-1}\times\C^{\varepsilon(k,m_1,...,m_t,s)},$$ 
when $(k-1)\nmid s$ and
$$U_{k,s,m_1,...,m_t}^{\l\neq1,c=0}=\C\setminus \{0,1\}\times(\C^*)^{t-1}\times\C^{\varepsilon(k,m_1,...,m_t,s)},$$
 $$U_{k,s,m_1,...,m_t}^{\l=1}=(\C^*)^{t-1}\times\C^{\varepsilon(k,m_1,...,m_t,s)}\times\C,$$
when $(k-1)\mid s$.
In case $(k-1)\mid s$ the spaces $U_{k,s,m_1,...,m_t}^{\l\neq1,c=0} $ and $U_{k,s,m_1,...,m_t}^{\l=1}$
appear as subspaces of
$$U_{k,s,m_1,...,m_t}=\C^*\times(\C^*)^{t-1}\times\C^{\varepsilon(k,m_1,...,m_t,s)}\times\C,$$
which parameterizes germs $  (z, \z) \mapsto (\l \z^s z +P(\z)+c\z^{\frac{sk}{k-1}}, \z^k)$
in modified normal form.
By the previous section all these germs have the same {\em logarithmic type}, 
i.e. the same directed dual graph of rational curves for the associated intermediate surfaces. 
Conversely, we have seen that every intermediate surface with such  configuration of curves corresponds to a germ of this type. 
Moreover, one may perform the blow-ups and the glueing
over the parameter space $U_{k,s,m_1,...,m_t}$ thus obtaining a family 
$\S_{k,s,m_1,...,m_t}\to U_{k,s,m_1,...,m_t}$ of intermediate surfaces over $U_{k,s,m_1,...,m_t}$.
It is clear that $\S_{k,s,m_1,...,m_t}$ is a complex manifold of dimension
$\dim U_{k,s,m_1,...,m_t}+2=t+\varepsilon(k,m_1,...,m_t,s)+\delta+2$, where $\delta=1$ if $(k-1)\mid s$ and otherwise $\delta=0$.
The projection $\S_{k,s,m_1,...,m_t}\to U_{k,s,m_1,...,m_t}$ is proper and smooth, since locally around each point of 
$\S_{k,s,m_1,...,m_t}$ it looks like the projection on $U_{k,s,m_1,...,m_t}$ of a small open subset of
$\hat B \times U_{k,s,m_1,...,m_t}$. 
See also our Appendix for a discussion on the family of the open surfaces 
$S\setminus D$.

Consider now an intermediate surface $S$ with maximal effective reduced divisor $D$. 
We would like to consider
families of logarithmic deformations as in \cite{Kaw78}. 
However the definition of that paper requires that $D$ be a 
subspace with {\em simple normal crossings}. 
But when the cycle of curves of $D$ is reduced to only one curve $C$, this
condition is not satisfied. 
In this case we blow up the singularity of $C$ on $S$ and work with the blown-up surface
instead. We may then apply Theorem 2 of \cite{Kaw78} to compare the logarithmic derformations. 
In the sequel we shall
work with the blown-up surface, without mentioning it explicitely again. In particular the new $D$ will have simple normal crossings. 
(Another way to avoid the singularity of $C$ would be to look at a non-ramified double cover of $S$ instead of $S$.)

A {\em family of logarithmic deformations} for the pair $(S,D)$ is a 6-tuple \\
$(\S,\D,\pi,V,v,\psi)$, where $\D$ is a divisor on $\S$, $\pi:\S\to V$ is a proper smooth morphism of complex spaces, which is locally a projection as well as its restriction to $\D$, $v\in V$
 and $\psi: S\to\pi^{-1}(v)$ is an isomorphism restricting to an isomorphism
 $S\setminus D\to\pi^{-1}(v)\setminus \D$. Let $T_S(-\log D)$ be the {\em logarithmic tangent sheaf} of $(S,D)$.
It is the dual of the sheaf $\Omega_S(\log D)$ of logarithmic differential 1-forms on $(S,D)$.
By \cite{Kaw78} versal logarithmic deformations of $(S,D)$ exist and their tangent space is $H^1(S,T_S(-\log D))$.
The space $H^2(S,T_S(-\log D))$ of obstructions vanishes by Theorem 1.3 of \cite{Nak90}. In particular the basis of the versal logarithmic deformation of a pair $(S,D)$ is smooth. On the other side 
 $H^0(S,T_S(-\log D))=H^0(S,T_S)$ and this space is at most one dimensional.

\begin{Th}\label{versal}
With the above notations we have:\\
\begin{itemize}
\item
 If $(k-1)$ does not divide $s$ the family \\
$$\S_{k,s,m_1,...,m_t}\to U_{k,s,m_1,...,m_t}$$
is logarithmically versal around every point of  
$U_{k,s,m_1,...,m_t}$. 
\item 
If $(k-1)\mid s$ the restriction of
the family $$\S_{k,s,m_1,...,m_t}\to U_{k,s,m_1,...,m_t}$$
to $U_{k,s,m_1,...,m_t}^{\l\neq1,c=0}$ is logarithmically versal
around every point of \\
$U_{k,s,m_1,...,m_t}^{\l\neq1,c=0}$. 
\item
If $(k-1)\mid s$ 
the family $\S_{k,s,m_1,...,m_t}\to U_{k,s,m_1,...,m_t}$
is logarithmically versal
around every point of $U_{k,s,m_1,...,m_t}^{\l=1}$.
\end{itemize}
\end{Th}
\begin{proof}
   We start with the case $(k-1)\nmid s$.
Take $(S,D)$ a pair as above and 
$ \varphi(z, \z):= (\l \z^{s} z +P(\z), \z^{k})$ 
an associated
 germ in  normal form (CG). We see $\varphi$ as a point in $U_{k,s,m_1,...,m_t}$.
Let $(\S',\D',\pi,V,v,\psi)$ be the logarithmically versal deformation of the pair $(S,D)$.
Then the family $\S_{k,s,m_1,...,m_t}\to U_{k,s,m_1,...,m_t}$ is obtained
from $\S'\to V$ by base change by means of a map 
$F:(U_{k,s,m_1,...,m_t},\varphi)\to(V,v)$.
By Theorem 4.\ref{favre} $F$ must have finite fibres. Thus
$\dim U_{k,s,m_1,...,m_t}\leq \dim V$. Since every point $v'$ of $V$ is covered by a similar map
$(U_{k,s,m_1,...,m_t},\varphi')\to(V,v')$ one gets $\dim U_{k,s,m_1,...,m_t}= \dim V$. 

Next we  show that $F$ is injective near $\varphi$ hence locally biholomorphic. Assume the contrary. 

Then by Theorem 4.\ref{favre}  there exists
a root of unity $\e$ of order $q$ with $1\neq q\mid (k-1)$ and a sequence
$(\varphi_n)_{n\in\N}$, $ \varphi_n(z, \z):= (\l_n \z^{s} z +P_n(\z), \z^{k})$, converging to $\varphi$ in $U_{k,s,m_1,...,m_t}$
such that for $ \tilde{\varphi}_n(z, \z):= (\e^s \l_n \z^{s} z +\e^{-m_1} P_n(\e\z), \z^{k})$ one has
$ \tilde{\varphi}_n\to\varphi$ and $F(\tilde{\varphi}_n)=F({\varphi}_n)$ for all $n\in\N$.
Let $r=-m_1- \lceil\frac{m_1}{q} \rceil q$ and $\chi(z,\z)=(\e^{-r}z,\e\z)$.
Then $\tilde{\varphi}_n=\chi^{-1}\circ {\varphi}_n\circ\chi$ for all $n\in\N$. 
Thus $\chi$ induces an automorphism 
of the germ $\varphi$ and in fact a non-trivial automorphism of the surface $S=S_{\varphi}$ as in the proof of
Remark 4.\ref{automorph}. On the other side $\chi$ induces for each $n\in\N$  the isomorphisms
$S_{\varphi_n}\cong S'_{F(\varphi_n)}\cong S_{\tilde{\varphi}_n}$ which by construction converge to the identity on 
$S=S_{\varphi}$. This is a contradiction.

The second assertion of the theorem has a completely analogous proof.

For the third  we have  $(k-1)\mid s$ and we consider the pair $(S,D)$ associated to a germ  
$  \varphi(z, \z) := (\l \z^s z +P(\z)+c\z^{\frac{sk}{k-1}}, \z^k)$ in normal form (CG),
 $(\S',\D',\pi,V,v,\psi)$ its logarithmically versal deformation and \\
$F:(U_{k,s,m_1,...,m_t},\varphi)\to(V,v)$ the induced morphism as before. From the second part of the theorem
it follows that\\
 $\dim U_{k,s,m_1,...,m_t}=\dim H^1(S,T_S(-\log D))=\dim V$. 
Since the fiber over $v$ is finite we only need to check that $F$ is injective.
If it was not, then
by the discussion of the first case ramification should occur along the divisor $U_{k,s,m_1,...,m_t}^{\l=1}$.
But this is not possible, because one can define the function $\l$ on the base of every family of
intermediate surfaces of this logarithmic type, in particular also on $V$. Indeed, for such a surface $S$ one 
defines $\l(S)$ as the unique twisting factor such that $H^0(S,T_S\otimes L_{\l(S)})\neq 0$, see for instance \cite{OTZ01}.
\end{proof}

We now fix a logarithmic type of intermediate surfaces, i.e. we fix the directed dual graph of the maximal reduced
 divisor of rational curves on such a surface, and want  to describe corresponding logarithmic moduli spaces.
By our previous discussions it is enough to fix one set $k,s,m_1,...,m_t$ of adapted numerical invariants and look 
at polynomial germs parameterized by $U_{k,s,m_1,...,m_t}$,
 $U_{k,s,m_1,...,m_t}^{\l\neq1,c=0}$, $U_{k,s,m_1,...,m_t}^{\l=1}$. 
Although several such sets of invariants may correspond to the desired moduli space, they all will be related by
the transformations $\tau_d$ described in sections \ref{highindex}, \ref{germs}.
 By Theorem 4.\ref{favre} we get a natural action of $\Z/(k-1)$ on $U_{k,s,m_1,...,m_t}$ which permutes conjugated
germs: through a generator of $\Z/(k-1)$ a germ 
$$ \varphi(z, \z) := (\l \z^s z +P(\z)+c\z^{\frac{sk}{k-1}}, \z^k)$$
is mapped to 
$$  \varphi(z, \z) := (\e^s \l_n \z^{s} z +\e^{-m_1} P_n(\e\z)+\e^{\frac{sk}{k-1}-m_1}cz^{\frac{sk}{k-1}}, \z^k),$$
where $\e$ is a primitive $(k-1)$-th root of the unity.

As a corollary to Theorem 7.\ref{versal} we get
\begin{Th}\label{moduli} 
Fix $k$, $s$ and a type $(m_1,...,m_t)$ for polynomial germs in normal form (CG).
Set $j=m_1$ as before.
\begin{itemize}
\item
 When $j<\max(s,k-1)$ the natural action of $\Z/(k-1)$ on  $$U_{k,s,m_1,...,m_t} ,
 U_{k,s,m_1,...,m_t}^{\l\neq1,c=0}\ \ {and} \ \ U_{k,s,m_1,...,m_t}^{\l=1}$$ is effective. 
\item 
In the remaining case, i.e. when $j=k-1=s$, the natural action of $\Z/(k-1)$ 
is effective on $U_{k,s,m_1,...,m_t}^{\l=1}=\C$ and trivial on 
$U_{k,s,m_1,...,m_t}^{\l\neq1,c=0}= \C\setminus\{ 0,1\}$.
\end{itemize}

The quotient spaces $U_{k,s,m_1,...,m_t}/(\Z/(k-1))$ when $(k-1)\nmid s$,
and \\
$U_{k,s,m_1,...,m_t}^{\l\neq1,c=0}/(\Z/(k-1))$, $U_{k,s,m_1,...,m_t}^{\l=1}/(\Z/(k-1))$,
when $(k-1)\mid s$, are coarse logarithmic moduli spaces for intermediate surfaces of the given logarithmic type
without, respectively with, non-trivial holomorphic vector fields.

These spaces are fine moduli spaces if and only if either the corresponding action of $\Z/(k-1)$
is trivial or this action is free.

The natural action of $\Z/(k-1)$ on either of the spaces $$U_{k,s,m_1,...,m_t},
U_{k,s,m_1,...,m_t}^{\l\neq1,c=0} \ \ {\rm and} \ \ U_{k,s,m_1,...,m_t}^{\l=1}$$ is free if and only if
$\gcd(k-1,s,m_2-j,...,m_t-j)=1$.
\end{Th}
\begin{proof}
The assertions on the effectivity of the action are immediately verified.

When the action is trivial it is clear that the corresponding families over
$U_{k,s,m_1,...,m_t}$,
 $U_{k,s,m_1,...,m_t}^{\l\neq1,c=0}$, $U_{k,s,m_1,...,m_t}^{\l=1}$ are universal.

When the action is free the families over
$U_{k,s,m_1,...,m_t}$,
 $U_{k,s,m_1,...,m_t}^{\l\neq1,c=0}$, $U_{k,s,m_1,...,m_t}^{\l=1}$
descend to the moduli spaces since we can extend the conjugation $\chi(z,\z)=(\e^{-r}z,\e\z)$
from the proof of Theorem 7.\ref{versal}
 to the whole family $\S_{k,s,m_1,...,m_t}\to U_{k,s,m_1,...,m_t}$.

When the action is not free but effective we see as in the proof of Theorem 7.\ref{versal}  that such a family
around a fixed point for some non-trivial subgroup of $\Z/(k-1)$ cannot descend to the quotient.

Suppose now that $d\mid\gcd(k-1,s,m_2-j,...m_t-j)$. It is easy to check that a germ of the form
$  \varphi(z, \z) := (\l \z^s z +\z^{m_1}+c_{m_2}\z^{m_2}+...+c_{m_t}\z^{m_t}, \z^k)$ is a fixed point for the
action of the subgroup of order $d$ of $\Z/(k-1)$.

Conversely, the existence of a fixed point for the
action of the subgroup of order $d$ of $\Z/(k-1)$ implies  $d\mid\gcd(k-1,s,m_2-j,...m_t-j)$.
\end{proof}

We have treated surfaces with and surfaces without vector fields separately in order to avoid non-separation phenomena. If we look for example at the families
$$\varphi_{1,\l}(z, \z) := (\l \z^{k-1} z +\z+c_1\z^k, \z^k),$$
$$\varphi_{2,\l}(z, \z) := (\l \z^{k-1} z +\z+c_2\z^k, \z^k),$$
for $c_1\neq c_2$ fixed and $\l$ varying in $\C^*$,
we see that $\varphi_{1,\l}$, $\varphi_{2,\l}$ are conjugated to one another if and only if 
$\l\neq 1$.

\section{Appendix: Deformation families for $\widetilde {S \setminus D}$ }

In this Appendix we shall construct deformation families for the universal covers of the open surfaces
$S\setminus D$ where $S$ is an intermediate surface and $D$ its maximal, effective, reduced divisor.
As a side product of our investigations we obtain examples of fibered analytic spaces together with groups
acting
holomorphically on them, such that the actions are free and properly discontinuous when restricted to each fiber
but not on the total spaces.

For simplicity we shall restrict ourselves to intermediate surfaces of index one.

Let $S$ be an intermediate surface and $D$ its maximal, effective, reduced divisor. 
We shall relax the conditions on the normal form (FG2) for the generators of the fundamental group of
$S\setminus D$ to get the following normal form:
$$\left\{
\begin{array}{lcl}
g_{\gamma}(z,w)&=&(z,w+ 2 \pi i)\\
g(z,w)&=&(\lambda z+a_0+Q(e^{-w}),kw),
\end{array}
\right. \leqno{(FG4)}
$$
where $k\geq 2$, $l\geq 1$, $Q=Q(\zeta):= \sum_{m=l}^{lk-1} b_m \zeta^m$ is a  polynomial
such that $(b_{lk-k+1},...,b_{lk-1})\neq 0$, 
$\gcd \{k,m \mid b_m \neq 0\} =1$. As before $(\lambda-1)a_0=0$. 
The difference from (FG2) is that we don't fix the leading coefficient $b_{\sigma}$ of $Q$ to be $1$
and require only that $lk-k<\sigma<lk$. In fact 
a conjugation by $(z,w)\mapsto(cz,w)$, with $c\in\C^*$, leads us to the generators
$$\left\{
\begin{array}{lcl}
g_{\gamma}(z,w)&=&(z,w+ 2 \pi i)\\
g(z,w)&=&(\lambda z+ca_0+cQ(e^{-w}),kw),
\end{array}
\right. .
$$
In the following we shall fix $k$ and $l$ and move the parameters \\
$(\l,a_0,b_l,...,b_{lk-1})\in \A\times\B$ of the normal form
(FG4). 
Here $\A:=\{(\l,a_0)\in\C^*\times\C \ \mid \ (\lambda-1)a_0=0\}$,
$\B:=\{ (b_l,...,b_{lk-1})\in\C^{lk-l} \mid \ (b_{lk-k+1},...,b_{lk-1})\neq 0\}$. 
For each fixed point $(\l,a_0,b)\in\A\times\B$
the operation of the group $\mathbb Z \ltimes \mathbb Z[1/k]$  generated by 
$g_{\gamma}$ and $g$ in the normal form (FG4) on $\mathbb C\times  \mathbb H_{l} $ 
is free and properly discontinuous and the quotient is a surface $S_{\l,a_0,b}\setminus D_{\l,a_0,b}$.

We can define the type for (FG4) in the same way as for the normal form (FG2), see Definition 3.\ref{type},
but note
 that in this case the type need not be constant on the $\B$-component.

\begin{Lem}\label{easy}
Consider the coefficients $b_l,...,b_{lk-1}$ of the polynomial $Q$ appearing in the normal form (FG4) of type
$(n_1,...,n_t)$ and set $n'_d:=\max \{n\in\N \ \mid \ b_n\neq0, \ \frac{k}{d}\nmid n\}$ for each 
 $d\in \{ 1,...,k-1 \} $ with $d\mid k$. Then $\{ n_1,...,n_t\}=\{ n'_d \ \mid \ 1\leq d<k, \ d\mid k\}$.
\end{Lem}
\begin{proof}
Set $j_0=1$. Then $n_i=n'_{\frac{k}{j_{i-1}}}$ for all $1\leq i\leq t$. Hence
$\{ n_1,...,n_t\}\subset\{ n'_d \ \mid \ 1\leq d<k, \ d\mid k\}$.

Conversely, suppose  $\{ n'_d \ \mid \ 1\leq d<k, \ d\mid k\}\setminus \{ n_1,...,n_t\}\neq\emptyset$
and let $n'_{\tau}=\min\{ n'_d \ \mid \ 1\leq d<k, \ d\mid k\}\setminus \{ n_1,...,n_t\}$,
$i=\max\{ j \in\N \ \mid \ n_j> n'_{\tau}\}$.

Suppose first that $i<t$. Then $n_i>n'_{\tau}>n_{i+1}$ and $j_i\mid n'_{\tau} $ by the definition 
of $n_{i+1}$. On the other side $\frac{k}{\tau}\mid n_i$ and thus $\frac{k}{\tau}\mid j_i$
by the definition of
of $n'_{\tau}$. Hence $\frac{k}{\tau}\mid j_i\mid n'_{\tau}$: a contradiction.

When $i=t$, we get $\frac{k}{\tau}\mid n_j$ for all $1\leq j\leq t$, which implies
$\frac{k}{\tau}\mid \gcd (n_1,...,n_t)=1$: again a contradiction.
\end{proof}

\begin{Th} \label{properdiscont} There is an action of $\mathbb Z \ltimes \mathbb Z[1/k]$  generated by
the holomorphic automorphisms 
$g_{\gamma}$ and $g$ in the normal form (FG4) on $\mathbb C\times  \mathbb H_{l} \times\A\times\B$ which is compatible
 with the projection  $\mathbb C\times  \mathbb H_{l} \times\A\times\B\to \A\times\B $.
Let $T\subset\A\times\B$ be a connected analytic subspace. The restriction of the above action 
to $\mathbb C\times  \mathbb H_{l} \times T$ is properly discontinuous if and  only if  the type is constant on $T$.
In this case the quotient space fibers smoothly over $T$.
\end{Th}
\begin{proof}
We consider the subgroup $\G<\mathbb Z \ltimes \mathbb Z[1/k]$, 
$$\G:=\{g_{r,m}:=g^{-m}\circ g_{\gamma}^r \circ g^m \ \mid \ m\in\Z, \ r\in \Z\}$$
We have $\G\cong  \mathbb Z[1/k]$. Actually 
$$g_{r,m}(z,w)=\Big(z+ \sum_{j=0}^{m-1}\l^{-j-1}(\sum_{n=l}^{lk-1}b_ne^{-nk^jw}(1-\exp(\frac{2\pi irn}{k^{m-j}}))),
w+\frac{2\pi ir}{k^{m}}\Big)$$
for each point $(z,w)\in \mathbb C\times  \mathbb H_{l} $. 
Note that in the above formula for $g_{r,m}$ we may always reduce ourselves to the situation when $k\nmid r$.

One can show as in \cite{DOT01} p. 659, that the action of $\mathbb Z \ltimes \mathbb Z[1/k]$ on 
$\mathbb C\times  \mathbb H_{l} \times T$ is properly discontinuous if and only  the induced action of 
the subgroup $\G$
is properly discontinuous.

Suppose now that the type is constant on $T$, 
that the action of $\G$ on $\mathbb C\times  \mathbb H_{l} \times T$ is not properly discontinuous and let
$(g_{r_{\nu},m_{\nu}})_{\nu\in\N}$ be a sequence in  $\G$ which contradicts the proper discontinuity of the action.
Then $(\frac{2\pi ir_{\nu}}{k^{m_{\nu}}})_{\nu}$ will be a bounded sequence, so by passing to some subsequence we may
assume that $(m_{\nu})_{\nu}$ is a strictly increasing sequence. Again by passing to some subsequence if necessary
we may assume that $\gcd (k, r_{\nu})=:d$ is constant for $\nu\in\N$. 
Let 
$n'_d:=\max \{n\in\N \ \vert \ b_n\neq 0, \ \frac{k}{d}\nmid m\}$. Since the type is constant on $T$ and by the
 above lemma $n'_d$ is well-defined. 
Since $(g_{r_{\nu},m_{\nu}})_{\nu\in\N}$ contradicts the proper discontinuity
of the action, for $(\l,a_0)$ bounded on the $\A$-component of $T$ and $w$ bounded in $ \mathbb H_{l}$ the quantity
$$\sum_{j=0}^{m_{\nu}-1}\l^{-j-1}(\sum_{n=l}^{lk-1}b_ne^{-nk^jw}(1-\exp(\frac{2\pi ir_{\nu}n}{k^{m_{\nu}-j}})))$$
must  be equally bounded independently of ${\nu}\in\N$. 
It suffices now to remark that the term 
$$\l^{-m_{\nu}}b_{n'_d}e^{-n'_dk^{m_{\nu}-1}w}(1-\exp(\frac{2\pi ir_{\nu}n'_d}{k}))$$
is dominant in this expression. 
Indeed the factor
$\vert 1-\exp(\frac{2\pi ir_{\nu}n'_d}{k})\vert$ is bounded from below by 
$\vert 1-\exp(\frac{2\pi i}{k})\vert$.
 On the other side the exponent $n'_dk^{m_{\nu}-1}$ is the highest
appearing in a non-vanishing term of this sum, since $kn'_d\geq kl>kl-1$ and thus
$b_n=0$ for $n>kn'_d$. 
Hence the term $\l^{-m_{\nu}}b_{n'_d}e^{-n'_dk^{m_{\nu}-1}w}(1-\exp(\frac{2\pi ir_{\nu}n'_d}{k}))$
of the above sum goes to infinity more rapidly than the rest and the sum cannot be bounded: a contradiction.

Conversely suppose that the type is not constant on $T$ and consider the first $n_i$ appearing 
in the decreasing sequence $n_1$, $n_2$, ... of type, such that $b_{n_i}$ takes both zero and non-zero values on $T$.
It is clear that one can find an analytic arc $\gamma:\Delta\to T$ such that $b_{n_i}$ vanishes at $\gamma(0)$ but such that the type is constant on $\gamma(\Delta\setminus\{0\})$ and contains $n_i$. 
Let $d$ be such that $n_i=n'_d$ for the generic type as in Lemma 8.\ref{easy}.
For the non-generic type the new
$n'_d$, call it $n''_d$ will assume a strictly lower value, by definition.
Take $w=-1$, $z$ arbitrary, $r_{\nu}=d$ for all ${\nu}\in\N$ and $m_{\nu}={\nu}$.
By the previous argument, for ${\nu}>>0$, the dominant term at $\gamma(0)$
$$\l^{-{\nu}}b_{n''_d}e^{n''_dk^{{\nu}-1}}(1-\exp(\frac{2\pi id_{\nu}n''_d}{k}))$$
will be lower than
$$\l^{-{\nu}}b_{n'_d}e^{n'_dk^{{\nu}-1}}(1-\exp(\frac{2\pi id_{\nu}n'_d}{k}))$$
in case $\mid b_{n'_d} \mid$ is uniformly bounded from below by some positive constant.
We can now choose a sequence converging to $\gamma(0)$ in $\gamma(\Delta)$ in such a way that
the corresponding $b_{n'_d}(\nu) $ converge to zero at such a rate that the whole sum
$$\sigma(\nu):=\sum_{j=0}^{{\nu}-1}\l^{-j-1}(\sum_{n=l}^{lk-1}b_n(\nu)e^{nk^j}(1-\exp(\frac{2\pi idn}{k^{{\nu}-j}})))$$
converges to zero.
In order to see this we solve the equation
$$b_{n'_d}+\frac{\sum_{j=0}^{{\nu}-1}\l^{-j-1}(\sum_{n=l, n\neq n'_d}^{lk-1}b_n e^{nk^j}(1-\exp(\frac{2\pi idn}{k^{{\nu}-j}})))}{\sum_{j=0}^{{\nu}-1}\l^{-j-1} e^{n'_d k^j}(1-\exp(\frac{2\pi idn'_d}{k^{{\nu}-j}}))}=0
\leqno{(5)}
$$  
in $b_{n'_d}$ on $\gamma(\Delta)$ for each $\nu>>0$.
This is certainly possible if the $b_n$-functions are supposed to be constant for $n\neq n'_d$.
The solutions  $(b_{n'_d}(\nu))_{\nu}$ form in this case the desired sequence and the sums $\sigma(\nu)$ vanish.
In general we rewrite the equation (5) as 

$$b_{n'_d}+\sum_{n=l}^{ n'_d-1}C_n(\nu)b_n =0.
\leqno{(6)}
$$  
It is easy to see that the coefficients $C_n(\nu)$ converge to zero as $\nu$ tends to infinity and thus (6) will have a solution on $\gamma(\Delta)$ for $\nu>>0$. This closes the proof.
\end{proof}


\end{document}